\documentclass{article}
\usepackage{irule,amsmath,amssymb}

\newcommand{\gad}{\ga\To G }
\newcommand{\ip}{\mathrm{IPC}\e}
\newcommand{\dwn}{\downarrow}

\newcommand{\mcal}{\mathcal{M}}
\newcommand{\rf}[1]{(\ref{#1})}
\newcommand{\e}{\epsilon}

\renewcommand{\iff}{\leftrightarrow}

\def\To{\Rightarrow}
\def\al{\alpha}
\def\ga{\Gamma}
\def\de{\Delta}
\def\qed{\hskip2em$\dashv$}
\def\prf{{\bf Proof.} }
 \newtheorem{Def} {Definition}
\newtheorem{Cor} {Corollary}
 \newtheorem{Lem} {Lemma}
 \newtheorem{Th} {Theorem}
 \newtheorem{example} {Example}
 
\def\beq{\begin{equation}}
\def\eeq{\end{equation}}
\def\ben{\begin{enumerate}}
\def\een{\end{enumerate}}

\newcommand{\barr}{\begin{array}}
\newcommand{\earr}{\end{array}}
\newcommand{\bcor}{\begin{Cor}}
\newcommand{\ecor}{\end{Cor}}
\newcommand{\bd}{\begin{Def}}
\newcommand{\ed}{\end{Def}}
\newcommand{\bl}{\begin{Lem}}
\newcommand{\el}{\end{Lem}}
\newcommand{\bt}{\begin{Th}}
\newcommand{\et}{\end{Th}}
\newcommand{\bexa}{\begin{example}}
\newcommand{\eexa}{\end{example}}
\newcommand{\skep}{\hskip2em}
\begin{document}
 \title{Intuitionistic Existential Instantiation and Epsilon Symbol}
 \date{\today}
 \author{Grigori Mints}
\maketitle 
\begin{abstract}
A natural deduction system for intuitionistic predicate logic with existential instantiation rule presented
here uses Hilbert's $\e$-symbol. It is conservative over intuitionistic
predicate logic. We provide a completeness proof for a suitable Kripke
semantics, sketch an approach to a normalization proof, survey related
work and state some open problems. Our system extends intuitionistic
systems with $\e$-symbol due to A. Dragalin and Sh. Maehara. 
\end{abstract}
\section{Introduction}
In natural deduction formulations of classical and intuitionistic logic
existence-elimination rule is usually taken in the form
\[
\infer[\exists^-]{C}
 {\exists xA(x) & \deduce[\vdots]{C}{A(a)}
 }
\]
where $a$ is a fresh variable
Existential instantiation is a rule
\[
\infer[\exists i]{A(a)}{\exists xA(x)}
\]
where $a$ is a fresh constant. It is sound and complete (with suitable
restrictions) in the role of
existence-elimination rule in classical
predicate logic but is not sound intuitionistically, since it makes
possible for example the following derivation:
\[
\infer{\To(C\to\exists xA(x))\to\exists x(C\to A(x))}
 {\infer{C\to\exists xA(x)\To\exists x(C\to A(x))}
  {\infer{C\to\exists xA(x)\To C\to A(a)}
     {\infer[\exists i]{C\to\exists xA(x), C\To A(a)}
      {C\to\exists xA(x), C\To \exists xA(x)}
     }
  }
 }
\] 
There are several approaches in the literature to introduction of
 restrictions making this rule conservative over
intuitionistic predicate calculus.

We present  an approach using intuitionistic version of Hilbert's
epsilon-symbol and strengthening works  by A. Dragalin
\cite{dragalin} and  Sh. Maehara \cite{maehara} where $\e$-terms are
treated as partially defined.  Then 
a survey of extensions and related approaches including important paper by
K. Shirai \cite{shirai} is given and some problems
are stated. 

We do not include equality since in this
case adding of $\e$-symbol with natural axioms is not conservative
over intutionistic logic (\cite{mintsskolemDAN,osswald}). A simple
counterexample due (in other terms) to C. Smorynsky \cite{smorynski} is
\[
\forall x\exists yP(x,y)\to 
\forall xx'\exists yy'( Pxy \& Px'y'\& (x=x'\to y=y'))
\]
In our natural deduction system NJ$\e$
axioms and propositional inference rules are the same as in ordinary
intuitionistic natural deduction, the same holds for
$\forall$-introduction. The remaining  rules are as follows:
\beq
\infer[\exists i]{\ga\To F(\e xF(x))}{\ga\To\exists xF(x)}
\eeq
existential instantiation,
\beq
\label{enegrules2}
\infer{\ga,\de\To F(t)}
 {\ga\To t\dwn & \de\To \forall zF(z)}
\skep
\infer{\ga,\de\To \exists zF(z)}
{\ga\To t\dwn & \de\to F(t)}
\eeq
where 
\beq
\label{dwn}
\e xA(x)\dwn:=  \exists y(\exists xA(x)\to A(y)),
\eeq
and $t\dwn:=\top$ (the constant ``true'') if $t$ is a variable or
constant. 

Two semantics are given for  NJ$\e$, or more precisely to an equivalent
 Gentzen-style system $\ip$ (Section \ref{sgentzen}). The first
semantics, which is incomplete but convenient for a proof of conservative
extension property over IPC is defined in Section \ref{ssemantics1}.

The second semantics with a completeness proof for $\ip$ is given in
Section \ref{scompleteness}.

Section \ref{scutelim} presents a sketch of a possible proof of a normal
form theorem. Section \ref{sprevious} surveys some of the previous work
and Section \ref{sfurther} outlines some open problems. 

\section{Gentzen-style system IPC$\e$}
\label{sgentzen}
Let us state our  Gentzen-style rules for the intuitionistic predicate
calculus IPC$\e$ with $\e$-symbol. For simplicity we assume that the language does
not have function symbols except constants. Formulas and terms are defined
by familiar inductive definition plus additional clause:

If $A(x)$ is a formula then $\e xA(x)$ is a term. 

Derivable objects of IPC$\e$ are {\em sequents} $\ga\To A$ where $\ga$ is a
finite set of formulas, $A$ is a formula. This means in particular that
structural rules are implicitly included below. 

First, let's list the rules of the intuitionistic predicate calculus IPC
 without $\e$-symbol. 

 Axioms:
\[ \ga,A \Rightarrow A,\skep \ga,\bot\To  A\ .
\]
 Inference rules: 
\[
\infer[\To\&]{\ga\To A \& B }
 {\ga\To A  & \ga\To B }
\skep
\infer[\&\To] {A\&  B ,\gad} {A, B ,\gad}  
\]
\[\infer[\vee\To]{A\vee B ,\gad}{A,\gad &  B ,\gad} 
\skep
\infer[\To\vee]{\ga\To A\vee B } {\ga\To A } 
\ 
\infer{\ga\To A\vee B } {\ga\To B } 
\]  
\[
\infer[\to\To] {A\to B ,\gad}{\ga\To A &  B ,\gad}
\skep 
\infer[\To\to]{\ga\To A\to B }  {A,\Gamma\Rightarrow B } 
\]
\[
\infer[\To\exists]{\ga\To \exists xA(x)}{\ga\To A (t)}
\skep
\infer[\To\forall]{\ga\To \forall xA(x)}{\ga\To A (b)}
\]
\[
\infer[\exists\To]{\exists xA(x),\ga\To  G }{A(b),\ga\to G }
\skep
\infer[\forall\To]{\forall xA(x),\ga\To  G }{A(t),\ga\to G }
\]
\[
\infer[Cut]{\ga\to G }{\ga\to C & C,\ga\To G }
\]
For IPC$\e$ quantifier-inferences $\To\exists,\forall\To$ are modified by requirement that the term $t$
substituted in the rule should be ``defined'' (cf. \rf{dwn}).
\beq
\label{enegrules}
\infer[\forall\To]{\forall zF(z),\ga,\de\To \de}
 {\ga\To t\dwn & F(t),\de\To \de}
\skep
\infer[\To\exists]
      {\ga\To\de, \exists zF(z)}
{\ga\To \de,t\dwn & \ga\to\de, F(t)}
\eeq
  $\exists\To$-rule is also  changed for IPC$\e$:
\beq
\label{existsant}
\infer[\exists_\e\To]{\ga,\exists xA(x)\To  G }
      {A(\e xA(x))\ga\To  G }
\eeq
A routine proof shows  that IPC$\e$ is equivalent 
to a Hilbert-style system obtained  by weakening
familiar axioms for quantifiers to 
\[
(\e Q1)\ t\dwn \& \forall xA(x)\to A(t)
\]
\[(\e Q2)\ t\dwn \& A(t)\to \exists xA(x)  \]
and adding the axiom
\[
\exists xA(x)\to A(\e xA(x))
\]
\subsection{Equivalence of $\ip$ and   NJ$\e$}
\label{subequiv}
Let us remind that in natural deduction a sequent 
\[
A_1,\ldots,A_n\To A
\]
is used to indicate that A is deducible from assumptions $A_1,\ldots,A_n$. 
\bt
A sequent is provable in  NJ$\e$ iff it is provable in IPC$\e$. 
\et

\prf The proof is routine: every rule of one of these  systems is directly
derivable in the other system. Let's show derivations of the rules
$\exists i$ and  $\exists\To$ from each other using abbreviation 
 $e:=\e xF(x)$.
\[
\infer{\ga\To F(e)}
   {\ga\to \exists xF(x) & 
     \infer[\exists\To]{\exists xF(x)\To F(e)}{F(e)\To F(e)}
   }
\skep
\infer{\exists xF(x),\ga\To G}
 {\infer[\exists i]{\exists xF(x)\To F(e)}{\exists xF(x) \To\exists xF(x)}
&
  \infer{\ga\To F(e)\to G}{F(e),\ga\To G}
 }
\]
\section{A Kripke Semantics for Intuitionistic $\e$-symbol}
\label{ssemantics1}
To prove that IPC$\e$ is conservative over IPC we present an incomplete 
semantics modifying a semantics  from \cite{dragalin}. The main
modification  is in the definition of $t\dwn$ and 
treatment of  atomic formulas containing $\e$-terms $\e x A$. 

\bd
Let $w$ be a world in a Kripke model. Denote 
\[
\e xA(x)\dwn w:\equiv\ w \models \e xA(x)\dwn . 
\]
We say that a term  $\e xA$ is {\em defined} in $w$  iff $\e xA(x)\dwn w$.
\ed
Symbol $\bot$ in  in  next definition  indicates  the condition
\rf{undefinedelement} below.
\bd
An intuitionistic Kripke $\e\bot$-model (or simply model in this section)
\[
\mcal=(W,<,D,\models,V)
\]
has to satisfy the following conditions:

 $(W,<)$ is a  Kripke frame with a strict partial ordering $<$,

 $D$ is a domain function assigning to every $w\in W$ a  non-empty set $D(w)$
monotone with respect to $<$,

 $w\models A$ is a relation between  worlds $w\in W$ and  atomic formulas
$A$ with constants from 
\[
D:= \cup_{w\in W}D(w)
\]
 monotonic with respect to $\leq$ and such that 
\beq
\label{undefinedelement} 
w\not\models A \text{ if $A$ contains at least one constant in $D-D(w)$}.
\eeq
 $V$ is a valuation function assigning a constant 
$V(e,w)\in D$ to any $\e$-term $e$ (possibly containing constants from
 $D$) and $w\in W$. 

The relation $\models$ is extended to composite formulas in the familiar
way. The components of an $\e$-model have to satisfy following
conditions. 
\beq
\label{eeinside}
V(\e xB(x,\e y C),w)=V(\e xB(x,V(\e yC,w)),w),
\eeq
\beq
\label{eainside}
w\models A(\e y C)\iff w\models A(x,V(\e yC)),
\eeq
where substitution of $\e yC$ is safe, that is no free variable of 
$\e yC$ becomes bound. Also 

if $e\dwn w$ for a term $e:=\ \e xA(x)$, then 
\[V(e,w)\in D(w) \text{ and } V(e,w')=V(e,w) \text{ for every }
w'\geq w. 
\]
\ed
Note once more that an  atomic formula  $P(d_1,\ldots,d_n)$ is {\em false}
in a world $w$ 
if at least one of $d_i$ is not in $D(w)$.

This leads to incompleteness, for example formula 
\[
P(\e xP(x))\to \exists xP(x)
\]
is valid: if $\e xP(x)$ is undefined in a world $w$ then the premise is
false in $w$, otherwise the conclusion is true. However this formula is not
derivable, since it  implies $(C\to\exists xP(x))\to \exists x(C\to
P(x))$.  

The proofs of the next lemmata are routine.
\bl 
\label{lmonot} 
Let $t$ be a closed term, $A$ a closed formula with
constants from $D$. Then
\[ w\leq w' \to\ (t\dwn w\to t\dwn v' \&\ 
 (w\models A \to w'\models A))
\]
\el
\prf Simultaneous induction on $t,A$.
\bl
\label{lsoundness}
If $\ga$ is a set of formulas, $G$ a formula then $\ga\vdash G$ in $\ip$ implies 
$\ga\models G$.
\el
\prf Induction on derivations. Checking the rule $\exists\To$ uses the
fact that $\exists xA(x)$ implies $\e xA(x)\dwn$. It may be
interesting to check whether any other properties of the formula $t\dwn$ are
used. \qed
\bt
\label{tconservative}
If $A,B$ formulas without $\e$-symbol then $A\vdash B$ in $\ip$ implies
$A\vdash B$ in intuitionistic predicate logic IPC.
\et
\prf We need to prove that for every Kripke model
\[ \mcal_0=(W,<,D,\models_0)
\]
 for intutionistic predicate logic refuting $A\to B$ there is an 
 $\ip$-model refuting $A\to B$. Before applying the construction from
 \cite{dragalin}, let  us recall a refinement of a completeness theorem
 for intuitionistic predicate logic IPC.
\bl
The following additional requirements to the definition of Kripke frame
$(W,<,D)$ for
IPC are still complete:
\ben
\item
\label{esearch1}
$W$ is a countable tree with a root $\bf0$ such that each $w\in W$ except
  $\bf0$ has unique  immediate $<$-predecessor and the number of predecessors
  of $w$ is finite.
\item
\label{esearch2}
domains D(w) are strictly increasing: if $w<w'$ then $D(w)$ is a proper
subset of $D(w')$.
\een
\el
\prf The requirement \ref{esearch1} is satisfied by the canonical proof
search tree for a given sequent, see for example \cite{mintsshort}. 
To satisfy the second requirement, note  that  an infinite branch of the
canonical proof search tree  does not have ``leaf worlds'': for every $w\in W$ there
exists a $w'>w$. Now take a fixed element $e\in D(w_0)$ and duplicate it
by a fresh element, say $e_w$ in every world $w$. More precisely for the
new domain function $D'$ define 
\[e_w\in D'(w)-D'(w^-),
\]
where $w^-$ is the immediate predecessor of $w$. Let's extend the relation
$\models$ by identifying $e_w$ and $e$, more precisely define for atomic
formulas $P(c_1,\ldots,c_n)$ with constants $c_i\in D'(w)$
\[
w\models P(c_1,\ldots,c_n)\ := w\models P(c_1^-,\ldots,c_n^-)
\]
where $c_i^-=e$, if $c_i=e_w$ and $c_i^-=c_i$ otherwise. 
It is easily proved by induction on formulas that this property extends to
all formulas:
\[
w\models A(c_1,\ldots,c_n)\text{ implies } w\models A(c_1^-,\ldots,c_n^-)
\]
so that the new model verifies (and refutes) the same formulas. \qed

Proof of the Theorem \ref{tconservative}.
We extend the model for IPC satisfying the previous Lemma by the
definition of values for $\e$-terms without changing domains $D(w)$,  which is done by induction on
construction of the term.  Assume that the elements of $D$ are
well-ordered by a relation $\prec$ in some arbitrary way. In view of the
condition \rf{eeinside} it enough 
to define $V(\e xA,w)$ when $\e xA$ does not have proper non-closed
$\e$-subterms. In that case, 

if $\e xA(x)\dwn w$, take the  $<$-minimal element $v\leq w$ such that 
 $\e xA\dwn v$, then define 
\[V(\e xA(x),w):=\ \text{the $\prec$-first } d\in D(v) 
(v\models (\exists xA(x)\to A(d)))
\]
If not $\e xA(x)\dwn w$, define $V(\e xA(x),w)$ as the $\prec$-first
$d\in D-D(w)$.
\qed
\section{Completeness proof for $\ip$}
\label{scompleteness}
We prove that removing condition \rf{undefinedelement} but preserving
familiar monotonicity requirement 
\beq
\label{monformulas}
w\leq w'\to\ (w\models A\to w'\models A)
\eeq
leads to a complete semantics for $\ip$. 

For simplicity consider term models where individual domain $D(w)$ for
every world $w$ consists of terms, and evaluation function for terms is
identity: value of a term $t$ is $t$. In particular 
the value of $\e xA$ is $\e xA$. 

\bd
An intuitionistic Kripke (term) $\e$-model (or simply $\e$-model)
\[
\mcal=(W,<,D,\models,V)
\]
has to satisfy the following conditions.

 $(W,<)$ is a  Kripke frame with a strict partial ordering $<$,

 $D$ is a domain function assigning to every $w\in W$ a  non-empty set
$D(w)$ (of terms) 
monotone with respect to $<$, 

 $w\models A$ is a relation between  worlds $w$ and  atomic formulas
$A$ with constants from 
\[
D:= \cup_{w\in W}D(w)
\]
 monotonic with respect to $\leq$ . 

 $V$ is a valuation function assigning a constant 
$V(e,w)\in D$ to any $\e$-term $e$ (possibly containing constants from
 $D$) and $w\in W$. (In a term model $V(e,w)=e$).

The relation $\models$ is extended to composite formulas in the familiar
way. The components of an $\e$-model have to satisfy following
conditions. 
\beq
\label{eeinside1}
V(\e xB(x,\e y C),w)=V(\e xB(x,V(\e yC,w)),w)
\eeq
\beq
\label{eainside1}
w\models A(\e y C)\iff w\models A(x,V(\e yC))
\eeq
 where substitution of $\e yC$ is safe, that is no free variable of 
$\e yC$ becomes bound. Also 
if $e\dwn w$ for a term $e:=\ \e xA(x)$, then
\[V(e,w)\in D(w) \text{ and } V(e,w')=V(e,w) \text{ for every }
w'\geq w. 
\]
\ed
Let's present a completeness proof along familiar lines.
\bd
An infinite sequent  is a pair of sets $\ga,\de$ of formulas such that
there is an infinite number of variables not in $\ga\cup\de$. An infinite
sequent $w$ is written as $\ga\To\de$ and notation 
\[
w_a:=\ga,\ w_s:=\de
\]
is used for its antecedent and succedent. 

$L_w$ denotes the set of all terms and formulas with free variables and
constants occurring in formulas of $w$.

$D(w)$ is the set of all terms  $t\in L_w$ such that $(t\dwn)\in w_a$. In other
worlds $D(w)$ consists of all free variables and constants in $w$ plus all
$\e$-terms $\e xA(x)$ such that $\exists y(\exists xA(x)\to A(y))\in\ w_a$.

An infinite sequent $w$ is {\em consistent}, if it is underivable, that is
if no finite sequent $\ga\To\de$ with $\ga\subset w_a,\ \de\subset w_s$ is
derivable in IPC$\e$.

A consistent  infinite sequent $w$ is {\em maximal consistent} if $w_a\cup
w_s$ is the whole set of formulas in $L_w$.
\ed
\bl
Every consistent infinite sequent $w_0$  can be extended to a maximal
consistent sequent.
\el
\prf  Enumerate all formulas containing only free variables and constants
in $L_{w_0}$,  then add them one by one to $w_a$ or $w_s$ preserving
consistency. At the $n$-th stage of this
process a sequent $w_n$, an extension of $w_0$ by a finite number of
formulas   is generated. 

It cannot happen that at some stage $n$ of this process a
formula $A$ fits none of $w^n_a,w^n_s$, i.e., both of 
\[ w^n_a\To w^n_s,A_n\skep A_n,w^n_a\To w^n_s
\]
are inconsistent, since in that case $ w^n_a\To w^n_s$ is inconsistent by a cut
rule. \qed

Important example. If $w$ is $\forall xP(x)\To P(\e xQ(x))$ with $P\neq Q$, and the first
``undecided'' formula  is $\exists y(\exists xQ(x)\to Q(y))$ then this formula is added to
the succedent, since adding it to the antecedent results in an inconsistent
sequent. 
\bl
\label{linvertible}
Every maximal consistent infinite sequent $w$ is closed under invertible
rules of multiple-succedent version of $\ip$, that is under all rules
except $\To\forall, \To\to$. More precisely 
\[(A\&B)\in w_a\text{ implies } A\in w_a\text{ and } B\in w_a,
\]
\[(A\to B)\in w_a \text{ implies } A\in w_s\text{ or } B\in w_a,
\]
\[(A\vee B)\in w_a \text{ implies } A\in w_a\text{ or } B\in w_a,
\]
\[(\forall xA(x))\in w_a  \text{ implies } 
(\forall t\in D(w))(A(t)\in w_a)
\]
\[(\exists xA(x))\in w_a \text{ implies } 
  A(\e xA(x))\in w_a
\]
\[(A\vee B)\in w_s\text{ implies } A\in w_a\text{ and } B\in w_a,
\]
\[(A\& B)\in w_s \text{ implies } A\in w_a\text{ or } B\in w_a,
\]
\[(\exists xA(x))\in w_s \text{ implies } 
(\forall t\in D(w))(A(t)\in w_s)
\]
\el
\prf Suppose $(A\&B)\in w_a$. If $A\not\in w_a$ then by maximality $A\in w_s$. Therefore
$w$ is inconsistent.

Suppose $\forall xA\in\ w_a$. If $A(t)\not\in w_a$ for some $t\in D(w)$
then by maximality  $A(t)\in w_s$. Therefore $\forall xA\To A(t)$ is
derived by one application of the $\forall\To$-rule, and hence 
$w$ is inconsistent. Note that additional premise $t\dwn$ of this rule is
available by  $t\in D(w)$.

Other cases are similar. \qed
\bd
For infinite sequents $w,w'$ define 
\[w<w'\text{ iff } w_a\subseteq w'_a \text{ and }
D(w)\subseteq D(w')
\]
\ed
\bl
\label{lnoninvertible}
The set of maximal consistent sequents is closed under non-invertible
rules $\To\to,\To\forall$. More precisely,

For every maximal consistent sequent $w$, if $(A\to B)\in w_s$ then there exists a
maximal consistent sequent $w'>w$ with 
$A\in w'_a,\ B\in w'_s$.

For every maximal consistent sequent $w$, if $\forall xA(x)\in w_s$ then there exists a
maximal consistent sequent $w'>w$ with 
$A(a)\in w'_s$  for some variable $a$, $a\in D(w')$.
\el
\prf  If $(A\to B)\in w_s$ then  the sequent $A,w_a\To B$ is
consistent, since otherwise one application of the rule $\To\to$ leads to 
 inconsistency of $w$. Now extend $A,w_a\To B$ to a complete consistent sequent.

 If $\forall xA(x) \in w_s$ then  the sequent $w_a\To A(a)$ for a fresh
 variable $a$  is
consistent, since otherwise one application of the rule $\To\forall$ leads to
 inconsistency of $w$. Now extend $A,w_a\To B$ to a complete consistent
 sequent.
\qed
\bd[Canonical model]
Consider the following model 
\[
M=(W,<,V,\models).
\]
$W$ is the set of all maximal complete sequents, $<,V$ are as above,

$w\models A$ iff $A\in w_a$ for atomic formulas $A$. 
\ed
This definition implies that  $w\not\models A$ for atomic $A\in w_s$, since
otherwise $w$ is inconsistent. 
\bl
The relation $\models$ for atomic formulas and  the function $D$ is
monotonic.  
\el
\prf Consider only $D(w)$. Let $w<w'$. All variables and constants in $D(w)$ are in $D(w')$ by
the definition of $<$. Assume $\e xA(x)\in D(w)$, that is 
$\e xA(x)\dwn\ \in w_a$. Then $\e x A(x)\dwn\in w'_a$ by $w<w'$, and
hence $\e xA(x)\in D(w')$. \qed
\bl
For every formula $A\in L_w$
\ben
\item $A\in w_a$ implies $w\models A$,

\item $A\in w_s$ implies $w\not\models A$,
\een
\el
\prf Induction on formulas using Lemmata
\ref{linvertible},\ref{lnoninvertible}. For example, if $A\& B\in w_a$ then 
$A,B\in w_a$, therefore $w\models A, w\models B$ by induction hypothesis,
and hence $w\models A\&B$.

If $\forall xA\in w_s$ then there exists $w'>w$ such that $A(a)\in w'_s$
for some variable $a\in D(w')$. Therefore $w'\not\models A(a)$ and hence 
$w\not\models\forall xA(x)$. \qed
\bt
The system IPC$\e$ is sound and complete.
\et
\prf Soundness is checked as before. For completeness take arbitrary
underivable formula $A$, then extend sequent $\To A$  to a maximal
consistent set $w$. By previous Lemma $w\not\models A$. \qed

%
\section{Cut-free Formulation, Normal Natural Deduction}
\label{scutelim}
 It is plausible that completeness proof for the rules with cut given in a 
 previous section \ref{scompleteness} can be modified to provide
 completeness of a  cut-free formulation. As our examples above show,
 complete cut-elimination is impossible.  One has to admit cuts for
 formulas of the form $\e xA(x)\dwn $ where $\e xA(x)$ occurs in the
 conclusion. The following proof where $e:=\e xP(x)$ is an example.
\[
\infer[cut]{\forall x\neg P(x),P(e)\To}
 { 
  \infer{\forall x\neg P(x),P(e)\To\exists y(\exists xP(x)\to P(y))}
  {\infer{\forall x\neg P(x),P(e)\To\exists xP(x)\to P(0)}
   {\infer{\forall x\neg P(x),P(e),\exists xP(x)\To P(0)}
    {\exists xP(x)\To e\dwn
      &
     \neg P(e),P(e),\exists xP(x)\To P(0)
    }
   }
  }
&
  \infer{\exists y(\exists xP(x)\to P(y)),\forall x\neg P(x),P(e)\To}
    {\exists y(\exists xP(x)\to P(y))\To e\dwn
      &
     \neg P(e),P(e)\To 
    }
 }
\]
 Let's outline a possible proof that this restriction is complete. 

First, the definition of the  canonical model should be modified along the
lines of a proof by M. Fitting \cite{fitting} (cf. also
\cite{mintsshort}). Our definition of a 
complete consistent sequent in the section \ref{scompleteness} requires
that such a sequent $w$  contains every formula of its language $L_w$ as a member of
its antecedent or succedent. This requirement is weakened as follows.

 For any formula $F$ in $L_w$ either $F\in w_a\cup w_b$ or there is a 
{\em clash}: both sequents 
\[ w_a\To w_b,F \text{ and } F,w_a\To w_b
\]
 are cut-free derivable. This should provide  completeness of a
 multiple-sequent  cut-free formulation.  Then equivalence to a cut-free
 one-succedent formulation should be proved by pruning and permutation of
 inferences as in \cite{mintsshort}. Finally cut-free one-succedent
 derivations are transformed into a normal natural deductions as in 
 \cite{mintslinear}.

\section{Comparison with Previous Work}
\label{sprevious}
\subsection{System $IPC\Omega\e$}

Let $\exists \e xA(x):= \exists xA(x)$.

A. Dragalin's system  $IPC\Omega\e$ from \cite{dragalin} for a given language $\Omega\e$ is
obtained  by weakening familiar axioms for quantifiers
\[
(\e Q1)\ \exists t\& \forall xA(x)\to A(t)
\]
\[(\e Q2)\ \exists t\& A(t)\to \exists xA(x)  \]
and adding the axiom
\[
\exists xA(x)\to A(\e xA(x))
\]
  A. Dragalin in \cite{dragalin} tried to avoid as much as
possible dealing with a  value of  an $\e$-term in a world $w$ where the
term is not defined.  Values (in a given world $w$) are assigned only to $\e$-terms defined
in $w$, and many intermediate results are proved only for the case when all
relevant $\e$-terms are defined. Nevertheless soundness is 
 established for all formulas, without any restrictions. As pointed
 earlier, this system is not complete.

In Section \ref{ssemantics1} we changed the  definition of a model from
\cite{dragalin} to a more uniform version: $\e$-term 
$e$ which is not defined at the world $w$ is assigned a value at $w$, but
this value does not belong to the individual domain $D(w)$. To make this
possible, the Kripke frame underlying the  model and the domain function
should satisfy additional conditions that still guarantee
completeness. 

Let us consider other systems in the literature.
\subsection{Systems with $\exists y(\exists xA(x)\to A(y)$ as Existence Condition} 
In  systems due to  to Sh. Maehara and K. Shirai \cite{maehara,shirai},
instead 
of using $\exists xA(x)$ as a discriminating criterion, a weaker formula
$\exists y(\exists xA(x)\to A(y))$ is employed. 
This still allows to anticipate a correct future value of the term
$\e xA(x)$ in a world $w$ even if $\exists xA(x)$ fails in $w$.

Sh. Maehara treats weaker language than ours:  $\e xA(x)$ is a
syntactically correct term only if it is closed. He proves (using partial 
cut-elimination and other syntactic transformations) conservativity
over IPC of the rules 
\[
\infer[\exists_\e]{\ga,\de\To G}
      {\ga,\exists xA(x) & A(\e xA(x)),\de\To G}
\]
\beq
\label{enegrules1}
 \infer{\forall zF(z),\ga,\de\To G}
 {\ga\To t\dwn & F(t),\de\To G}
\skep
\infer{\ga,\de\To \exists zF(z)}
{\ga\To t\dwn & \de\to F(t)}
\eeq
where
\beq
\label{edefined}
\e xA(x)\dwn:= \exists y(\exists xA(x)\to A(y));\ a\dwn:=\top
\eeq
Here $\top$ is the constant true, $a$ is an arbitrary variable.

Note that the first of these rules contains a hidden cut. This
conservativity result is used to establish a kind of completeness theorem
for IPC over a modification of Kripke semantics, although this
modification is not stated explicitly. More precisely, Sh. Maehara proves 
Kripke-style soundness and completeness result for the relation $A\in
\alpha$ between formulas 
$A$ and complete consistent (in his sense) subsets $\al$ of the set of formulas. Only his  condition for
$\forall$ is not standard:

$\forall xA(x)\in \al\iff\ (\exists B)(B\in\al \&\forall\beta\forall t 
[B\in\beta\to (t\in D_\beta\to A(t)\in\beta)])$

To establish this condition he uses admissibility of the following rule in
his system:
\[
\infer{\forall xA(x)}{\exists y(\exists x\neg A(x)\to\neg A(y))\to A(\e
  x\neg A(x))}
\]
This rules  approximates equivalence 
\[
\forall xA(x)\iff A(\e x\neg A(x))
\]
which is valid only classically. 

K. Shirai \cite{shirai}  removes restriction to closed $\e$-terms. He considers a
language with the existence predicate  denote by $D$. Instead of rules
used by Maehara he considers following axioms:
\beq
\label{Dsufficient}
D(t),\exists y(\exists xA(x,t)\to A(y,t)\To D(\e xA(x))
\eeq
\[
D(t),\exists xA(x,t)\To A(\e xA(x,t),t)
\]
 plus  standard modifications of quantifier rules for the system with
 existence predicate D. 

He  proves conservativity of his system over IPC by a combination of
a partial cut-elimination and Maehara's argument.

D. Leivant \cite{leivant} and V. Smirnov \cite{smirnov} define logical
systems with $\e$-symbol conservative over IPC by requiring that
assumptions discharged in   natural 
deduction rules contain no $\e$-symbol. These systems are probably much
weaker than $\ip$. The system introduced by the author in
\cite{mintsepsilon} is certainly weaker than $\ip$: a sequent containing
subterm $\e xA(x,y)$ with a bound variable $y$ is syntactically correct
only provided $\forall y\exists xA(x,y)$ is a member of the antecedent.
\section{Further Work}
\label{sfurther}
Complete proof of cut-elimination for IPC$\e$ and of normal form theorem
for NJ$\e$. 

Give a syntactic proof of cut-elimination for IPC$\e$ and of normalization
for NJ$\e$.

Provide a semantics for the systems by Sh. Maehara and K. Shirai
\cite{maehara,shirai} and find out whether these systems admit
cut-elimination. It seems that the system by Shirai provides the most
general formulation of the idea that $\e$-terms is partially defined in
some arbitrary way. The restriction $D(t)$ allowing to quantify over value
of $t$ can be arbitrary predicate with the only condition
\rf{Dsufficient}. 


\begin{thebibliography}{99}
\bibitem{dragalin}  Dragalin, A. Intuitionistic Logic and Hilbert's
  $\e$-symbol, (Russian) Istoriia i Metodologiia Estestvennykh Nauk,
  Moscow, MGU, 1974, s. 78-84, republished in: Albert Grigorevich
  Dragalin, Konstruktivnaia Teoriia   Dokazatelstv I Nestandartnyi
  Analiz, s. 255-263, Moscow, Editorial Publ. )
\bibitem{fitting} Fitting, M., Intuitionistic Logic, Model Theory And
  Forcing, Amsterdam, North-Holland, 1969 
  \bibitem{leivant} Leivant, D., Existential instantiation in a
    system of natural deduction for intuitionistic arithmetics,
  Technical  Report ZW 13/73, Stichtung Mathematisch Centrum, Amsterdam,
  1973 
\bibitem{maehara} Maehara, Sh., A General Theory of Completeness Proofs,
  Ann. Jap. Assoc. Phil. Sci., 1970, no. 3, p. 242-256
\bibitem{mintsskolemDAN}  Mints, G., The Skolem Method  in Intuitionistic
  Calculi. Proc. Inst. Steklov, 121, AMS, 1974, p. 73-109
\bibitem{mintsnonskolem} Mints, G., Skolem Method of Elimination of
  Positive Quantifiers in Sequential 
Calculi, Soviet Math. Dokl. 7, no.4, 1966, 861-864 
\bibitem{mintslinear} Mints, G., Linear Lambda-terms and Natural
  Deduction, Studia Logica, 60, 1998, p. 209-231 
\bibitem{mintsshort} Mints, G., A Short Introduction to Intuitionistic Logic ,
  Kluwer Academic/ Plenum Publishers, 2000
\bibitem{mintsepsilon} Mints, G., Heyting Predicate Calculus with Epsilon  
  Symbol (Russian), Zapiski Nauchnykh Seminarov Leningradskogo
  Otdeleniya Matematicheskogo Instituta im. V. A. Steklova AN SSSR,
  Vol. 40, pp. 110-118, 1974, English Translation in \cite{mintsselected} 
  p. 97-104
\bibitem{mintsselected}  Mints, G., Selected Papers in Proof
  Theory, Bibliopolis/North-Holland, 1992
\bibitem{osswald} Osswald, H., \"Uber  Skolemerweiterungen in der
  Intuitionistichen Logik mit Gleichheit, Lecture Notes in Mathematics,
  1975, Volume 500, 1975, p. 264-266, 
\bibitem{shirai} Shirai, K. Intuitionistic Predicate Calculus with
  $\e$-symbol, 1971, no.4, p. 49-67
\bibitem{smorynski} Smorynski C., On Axiomatizing
  Fragments, J. Symbolic Logic, v. 42, no.4, 1977, p. 530-544
\bibitem{smirnov} Smirnov, V., Theory of
  quantification and E-calculi. In:  Essays on mathematical and
  philosophical logic, Essays on Mathematical and Philosophical Logic:
  Proceedings of the Fourth ...
 By Jaakko Hintikka, Ilkka Niiniluoto, Esa Saarinen, Proc. 4th
 Scand. Logic Symp., Kluwer, 1979,  p. 41-49
\end{thebibliography}
\end{document}